\documentclass[a4paper, 12pt, fleqn]{article}
\usepackage[utf8]{inputenc} %paquete para usar teclado español
\usepackage[T1]{fontenc}%to be able to hyphenate accents
\usepackage{amsfonts} %\mathbb
\usepackage{amssymb}
\usepackage{amsmath}%binom
\usepackage{graphicx} %paquete para trabajar con imágenes
\usepackage{listings} %code C,Java,...
\usepackage{amsthm} %thm, proof % ---------- COMMENT for IEEE ------------------
\usepackage{color}
\usepackage[ruled,vlined]{algorithm2e}  % para el pseudo-codigo

\usepackage{float}%for true here \begin{figure}[H]
\usepackage{tikz}
\usepackage{todonotes}%\todo{text}
\usepackage[bookmarks=true]{hyperref} %hiper-referencias, claro!

\newcommand\titulo{Simple, Fast and Practicable Algorithms for Cholesky, LU and QR Decomposition Using Fast Rectangular Matrix Multiplication}
\newcommand\autor{Cristóbal Camarero}
\hypersetup{
    unicode=false,          % non-Latin characters in Acrobats bookmarks
    pdftoolbar=true,        % show Acrobats toolbar?
    pdfmenubar=true,        % show Acrobats menu?
    pdffitwindow=true,      % page fit to window when opened
    pdftitle={\titulo},     % title
    pdfauthor={\autor},     % author
    pdfsubject={},			% subject of the document
    pdfnewwindow=true,      % links in new window
    pdfkeywords={},			% list of keywords
	pdfpagemode=None,		% avoid auto open bookmarks
    colorlinks=false,       % false: boxed links; true: colored links
    linkcolor=red,          % color of internal links
    citecolor=green,        % color of links to bibliography
    filecolor=magenta,      % color of file links
    urlcolor=cyan           % color of external links
}


\begin{document}
\title{\titulo}
\author{Cristóbal Camarero\\
Department of Computer Science and Electronics\\
Universidad de Cantabria, UNICAN, Spain.}
\date{}
\maketitle%

\begin{abstract}This note presents fast Cholesky/LU/QR decomposition algorithms with $O(n^{2.529})$ time complexity when using the fastest known matrix multiplication. The algorithms have potential application, since a quickly made implementation using Strassen multiplication has lesser execution time than the employed by the GNU Scientific Library for the same task in at least a few examples.

The underlaying ideas are very simple. Despite this, I have been unable to find these methods in the literature. 
\end{abstract}

\section{Introduction}

Cholesky and LU decompositions are used everywhere to solve linear systems; and they have multiple applications beyond that.
The problem is given a positive-definite matrix $A$ to find two triangular matrices $L$ and $U$ such that $A=LU$, $L$ being lower triangular and $U$ being upper triangular. When $A$ is symmetric, it is hold that $U=L^T$ after normalization, and $A=LL^T$ is called a Cholesky decomposition.
In the QR decomposition the problems is to find orthogonal $Q$ and upper triangular $R$ such that $A=QR$.
The new presented algorithms reduce their time complexity from $O(n^3)$ to $O(n^{2.529})$. Although in \cite{Bunch} it was presented an algorithm for LU decomposition with lower complexity, it was not very practical because of the large constants involved. In contrast, the present algorithms are very practical, since they have no additional large constants respect to the typical Cholesky/LU/QR algorithms except the ones being present on the matrix multiplication algorithm of choice.
The same strategy could also be applied to tridiagonalization of symmetric matrices, but then we would just remade the algorithm in \cite{Bischof}.

%The speed of the algorithm depends on the complexity of rectangular matrix multiplication; for this kind of algorithms see \cite{LeGall2012,LeGall2018}. Their aim is to multiply a $n\times n^\alpha$ with a $n^\alpha\times n$ matrix in $O(n^{2+\epsilon})$ for any $\epsilon>0$ and the greatest possible $\alpha$, which has recently been increased to $\alpha=0.31389$ \cite{LeGall2018}.
The speed of the algorithms depends on the capacity to perform rectangular matrix multiplications. Matrix multiplication algorithms are usually studied for square matrices and they have time complexity of $O(n^\omega)$, where $n$ is the matrix side and $\omega$ depends on the algorithm: $\omega=3$ for classical multiplication, $\omega=2.8074$ for Strassen's algorithm and $\omega=2.3728639$ for the asymptotically best known matrix multiplication algorithm \cite{LeGall2014}. For the present algorithms for Cholesky and LU decompositions, the interest lies on multiplying a matrix with dimensions $n\times s$ with a $s\times n$ matrix, where $s$ is smaller than $n$. In \cite{Knight} it is shown that such multiplication can be performed with $O(s^{\omega-2}n^2)$ arithmetic operations using square-based algorithms, but also that algorithms operating directly on rectangles may outperform the square based ones.
In \cite{LeGall2012,LeGall2018} rectangular matrix multiplication is studied more thoroughly. Let $s=n^k$, the $\omega$ notation may be extended as $\omega(k)$ to indicate that the multiplication of a $n\times s$ by a $s\times n$ can be performed in $O(n^{\omega(k)})$. Furthermore, for small values of $k$, the optimal complexity of $O(n^2)$ is obtained, which is currently at $\alpha=0.31389$, and thus $\omega(\alpha)=2$. Nevertheless, algorithms with good $\omega(k)$ for every $k$ are found in \cite{LeGall2018}. We highlight the value $\omega(0.5286)=2.057085$, which we will use later.

\section{The algorithms}

\begin{algorithm} \SetLine
\SetKwFunction{Function}{Function}
\KwData{A square, symmetric, and positive-definite matrix $A$;\\
The size step $s$.}
\KwResult{The lower triangular matrix $L$ such that $A=LL^T$.}

$n$:= size of $A$ \;
$L$:= zero matrix of size $n$ (could be done inplace of $A$)\;
$z$:= 1 \;

\For{$c$ from 1 to $n$}
{
	\If{$c=z+s$}
	{
		$R$:= submatrix of $L$ by [$c$ to $n$] and [$z$ to $c-1$] \;
		$S$:= $RR^T$, using a fast rectangular matrix multiplication algorithm \;
		\For{$i$ from $c$ to $n$}
		{
			\For{$j$ from $c$ to $n$}
			{
				$A_{i,j}$:= $A_{i,j}-S_{i-c+1,j-c+1}$ \;
			}
		}
		$z$:= $c$ \;
	}
	$L_{c,c}$:= $A_{c,c}$ \;
	\For{$k$ from $z$ to $c-1$}
	{
		$L_{c,c}$:= $L_{c,c} - L_{c,k}^2$ \;
	}
	$L_{c,c}$:= $\sqrt{L_{c,c}}$ \;
	\For{$i$ from $c+1$ to $n$}
	{
		$L_{i,c}$:= $A_{i,c}$ \;
		\For{$k$ from $z$ to $c-1$}
		{
			$L_{i,c}$:= $L_{i,c} - L_{i,k}L_{c,k}$ \;
		}
		$L_{i,c}$:= $\frac{L_{i,c}}{L_{c,c}}$ \;
	}
}

\Return{$L$}

% \eIf {$\beta=(t+1)+(t+1)i$ or $|r| \leq t$} {Return $\theta = q
% \beta$ \;} {Compute:\par \hspace{2ex}$ Q = \{ (q + h) \beta
% \pmod{\alpha} \ | \ h \in \{0, \pm 1 \pm i \} \}$ \; Find $\theta$
% such that $ |\theta| = \min\{ D_{\alpha}(x, \gamma) \ | \ x \in Q
% \}$ \; Return $\theta$ \; }
% \caption{Decoding in $\mathbb{Z}[i]_{\alpha}$} \label{algo:Gauss}
\caption{A fast Cholesky decomposition algorithm}\label{algo:fastCholesky}
\end{algorithm}

\begin{algorithm} \SetLine
\SetKwFunction{Function}{Function}
\KwData{A square matrix $A$ with nonzero leading principal minors;\\
The size step $s$.}
\KwResult{The lower triangular matrix $L$ and upper unitriangular matrix $U$ such that $A=LU$.}

$n$:= size of $A$ \;
$L$:= zero matrix of size $n$ \;
$U$:= zero matrix of size $n$ \;
$z$:= 1 \;

\For{$c$ from 1 to $n$}
{
	\If{$c=z+s$}
	{
		$R_L$:= submatrix of $L$ by [$c$ to $n$] and [$z$ to $c-1$] \;
		$R_U$:= submatrix of $U$ by [$z$ to $c-1$] and [$c$ to $n$] \;
		$S$:= $R_LR_U$, using a fast rectangular matrix multiplication algorithm \;
		\For{$i$ from $c$ to $n$}
		{
			\For{$j$ from $c$ to $n$}
			{
				$A_{i,j}$:= $A_{i,j}-S_{i-c+1,j-c+1}$ \;
			}
		}
		$z$:= $c$ \;
	}
	\For{$i$ from $c$ to $n$}
	{
		$L_{i,c}$:= $A_{i,c}$ \;
		\For{$k$ from $z$ to $c-1$}
		{
			$L_{i,c}$:= $L_{i,c} - L_{i,k}U_{k,c}$ \;
		}
	}
	\For{$i$ from $c$ to $n$}
	{
		$U_{c,i}$:= $A_{c,i}$ \;
		\For{$k$ from $z$ to $c-1$}
		{
			$U_{c,i}$:= $U_{c,i} - L_{c,k}U_{k,i}$ \;
		}
		$U_{c,i}$:= $\frac{U_{c,i}}{L_{c,c}}$ \;
	}
}

\Return{$L,U$}

% \eIf {$\beta=(t+1)+(t+1)i$ or $|r| \leq t$} {Return $\theta = q
% \beta$ \;} {Compute:\par \hspace{2ex}$ Q = \{ (q + h) \beta
% \pmod{\alpha} \ | \ h \in \{0, \pm 1 \pm i \} \}$ \; Find $\theta$
% such that $ |\theta| = \min\{ D_{\alpha}(x, \gamma) \ | \ x \in Q
% \}$ \; Return $\theta$ \; }
% \caption{Decoding in $\mathbb{Z}[i]_{\alpha}$} \label{algo:Gauss}
\caption{A fast LU decomposition algorithm}\label{algo:fastLU}
\end{algorithm}

\begin{algorithm} \SetLine
\SetKwFunction{Function}{Function}
\KwData{A matrix $A$;\\
The size step $s$.}
\KwResult{An orthogonal matrix $Q$ and right upper triangular $R$ such that $A=QR$.}

$n$:= size of $A$ \;
$M$:= a copy of $A$ \;
$Q$:= zero matrix of size $n$ \;
$z$:= 1 \;

\For{$c$ from 1 to $n$}
{
	\If{$c=z+s$}
	{
		$B_Q$:= submatrix of $Q$ by [$1$ to $n$] and [$z$ to $c-1$] \;
		$B_M$:= submatrix of $M$ by [$1$ to $n$] and [$c$ to $n$] \;
		$C$:= $B_Q^TB_M$, using a fast rectangular matrix multiplication algorithm \;
		$S$:= $B_QC$, using a fast rectangular matrix multiplication algorithm \;
		\For{$i$ from $1$ to $n$}
		{
			\For{$j$ from $c$ to $n$}
			{
				$M_{i,j}$:= $M_{i,j}-S_{i,j-c+1}$ \;
			}
		}
		$z$:= $c$ \;
	}
	$v$:= column $c$ of $M$ \;
	\For{$j$ from $z$ to $c-1$}
	{
		$u$:= column $j$ of $Q$ \;
		$v$:= $v - u(u^Tv)$ \;
	}
	$v$:=$\frac{v}{||v||}$ \;
	\For{$i$ from $1$ to $n$}
	{
		$Q_{i,c}$:= $v_i$ \;
	}
}
$R$ := $Q^TA$, using a fast square matrix multiplication algorithm \;
\Return{$Q,R$}
\caption{A fast QR decomposition algorithm}\label{algo:fastQR}
\end{algorithm}

The new Cholesky decomposition algorithm is shown as Algorithm \ref{algo:fastCholesky}, the new LU decomposition algorithm as Algorithm \ref{algo:fastLU}, and the new QR decompositon algorithm as Algorithm \ref{algo:fastQR}. All indices start at 1 and end at the size of the matrix, usually $n$, and all ranges are inclusive. The notation $A_{i,j}$ indicates access to the $(i,j)$-entry of the matrix $A$. The three algorithms are analogous, in the sense they have the same kind of variation respect to the classical algorithms, so let us comment from just the Cholesky viewpoint. It combines the ideas of the recursive Cholesky algorithm, where the lower right submatrix is updated by $-vv^T$ in each call; and of the Cholesky--Crout algorithm, which computes a column using the previous computed columns.
In this new algorithm there is a parameter $s$ indicating the maximum number of columns to process before updating the submatrix. Thus, when $s=1$ it particularizes to the recursive algorithm and when $s=n$ it becomes Cholesky--Crout algorithm.
The correctness of the algorithm is trivially reduced to the correctness of these two algorithms it generalizes.

The updates of the lower right submatrix are made by substractions of a $RR^T$ block, where $R$ has (or can be extended to) size $n\times s$. When $s\leq n^{0.313}$, this $RR^T$ product can be computed in $O(n^2)$ \cite{LeGall2018} as said in the introduction, but other values of $s$ also have interest. There is a total of $n$ iterations. In each iteration there are $O(ns)$ arithmetic operations that always execute. Once every $s$ iterations the first conditional executes, which contains the update of the submatrix in $O(n^2)$ arithmetic operations preceded by the rectangular matrix multiplication.

Therefore, if $s=n^{0.313}$ then the total volumen of computation is of $O(n^2s+\frac{n}{s}n^2)=O(n^{2.687})$ arithmetic operations.
Using the fastest matrix multiplication algorithm for square matrices, $\omega=2.3728639$ \cite{LeGall2014}, we get complexity $O(n^2s+\frac{n}{s}s^{\omega-2}n^2)$. Its optimum is at $s=n^{0.61462815}$, where the algorithm becomes $O(n^{2.61462815})$.
But it is better to use the value $\omega(0.5286)=2.057085$, which leads to $s=n^{0.5285425}$ and time complexity of $O(n^2s+\frac{n}{s}n^{\omega(k)})=O(n^{2.5285425})$.
Further improvements to rectangular matrix multiplication could in theory reduce the final complexity of these decomposition algorithms to $O(n^{2.5})$.

In practice those multiplication algorithms with small complexity incur in huge constants. A good alternative is Strassen algorithm, which can be used in rectangular matrices \cite{Knight}. Taking $s=n^{0.8385}$, the whole algorithm becomes with Strassen multiplication $O(n^2s+\frac{n}{s}s^{\omega-2}n^2)=O(n^{2.8385})$, where $\omega=\log_2(7)=2.8074$ is the exponent in the cost of Strassen multiplication for square matrices.

An implementation of the algorithms using Strassen multiplication has proved to spend about 48\% less time than GNU Scientific Library (GSL) \cite{GSL} Cholesky implementation,  about 57\% less time than GSL LU implementation, and about 87\% less time for the QR.
For $n=2000$, taking $s=200$ and employing two Strassen recursions (and the remaining multiplications made with GSL) it took 1.14867 seconds, against the 2.28978 seconds of pure GSL. The test matrix was a symmetric one with random entries and large values in the diagonal.
For $n=4000$, we take $s=400$ and employ three Strassen recursions. This consumes 8.99871 seconds, while GSL consumes 17.0262 seconds.
The analogous for LU is, in $n=2000$, 1.51041 seconds vs 3.53489 seconds; and in $n=4000$, 12.203 seconds vs 28.7355 seconds; both with $s=200$. Although one must note that GSL LU also computes a permutation matrix.
The analogous for QR is, in $n=2000$, 7.12522 seconds vs 51.8309 seconds; and in $n=4000$, 64.3856 seconds vs 554.165 seconds; the first with $s=100$ and the second with $s=200$.

It is also remarkable that the $s$ parameter could be varied during the execution. Modifications of this kind could lead to some improvements in practice.

\section{Acknowledgments}
I would like to give thanks to Luis M. Pardo for reading a draft of this note and giving some advice.

The author is currently with a Juan de la Cierva--Formación contract (FJCI-2017-31643) by the Ministerio de Ciencia, Innovación y Universidades of Spain, which is cofinanced by the Department of Computer Science and Electronics of the Universidad de Cantabria in Spain.

%\begin{thebibliography}{99}
%\bibitem{ISIT10} C. Mart\'{\i}nez, C. Camarero, R. Beivide ``Perfect Graph Codes over Two Dimensional Lattices". Accepted for publication at 2010 IEEE Internations Symposium on Information Theory.
%\end{thebibliography}
\bibliographystyle{plain}
\bibliography{main}

\begin{thebibliography}{1}

\bibitem{Bischof}
Christian~H. Bischof, Bruno Lang, and Xiaobai Sun.
\newblock A framework for symmetric band reduction.
\newblock {\em ACM Trans. Math. Softw.}, 26(4):581--601, December 2000.

\bibitem{Bunch}
James~R. Bunch and John~E. Hopcroft.
\newblock Triangular factorization and inversion by fast matrix multiplication.
\newblock {\em Mathematics of Computation}, 28(125):231--236, 1974.

\bibitem{LeGall2012}
Fran\c{c}ois~Le Gall.
\newblock Faster algorithms for rectangular matrix multiplication.
\newblock In {\em 2012 IEEE 53rd Annual Symposium on Foundations of Computer
  Science}, pages 514--523, October 2012.

\bibitem{LeGall2018}
Fran\c{c}ois~Le Gall and Florent Urrutia.
\newblock {\em Improved Rectangular Matrix Multiplication using Powers of the
  {C}oppersmith-{W}inograd Tensor}, pages 1029--1046.

\bibitem{GSL}
Brian Gough.
\newblock {\em GNU Scientific Library Reference Manual - Third Edition}.
\newblock Network Theory Ltd., 3rd edition, 2009.

\bibitem{Knight}
Philip~A. Knight.
\newblock Fast rectangular matrix multiplication and {QR} decomposition.
\newblock {\em Linear Algebra and its Applications}, 221:69 -- 81, 1995.

\bibitem{LeGall2014}
Fran\c{c}ois Le~Gall.
\newblock Powers of tensors and fast matrix multiplication.
\newblock In {\em Proceedings of the 39th International Symposium on Symbolic
  and Algebraic Computation}, ISSAC '14, pages 296--303, New York, NY, USA,
  2014. ACM.

\end{thebibliography}

\end{document}